\documentclass[letter,12pt]{article}%
\usepackage{amsmath,bm}
\usepackage{amsfonts}
\usepackage{amssymb}
\usepackage{graphicx}
\usepackage{rotating}
\usepackage[height=9in,left=1in,right=0.75in,bottom=1in]{geometry}
\usepackage[breaklinks=true,bookmarksopen=true,colorlinks=true,citecolor=blue]%
{hyperref}
\usepackage[authoryear]{natbib}
\usepackage{color}
\usepackage{colortbl}
\usepackage{paralist}
\usepackage{amsmath}%
\providecommand{\U}[1]{\protect\rule{.1in}{.1in}}
\RequirePackage{hypernat}

\catcode`~=11
\newcommand{\urltilde}{\kern -.15em\lower .7ex\hbox{~}\kern .04em}
\catcode`~=13
\def \@seccntformat#1{\csname the#1\endcsname.\quad}
\makeatother
\numberwithin{equation}{section}
\hypersetup{pdftitle={Escanciano}, pdfsubject={Robust Two-Step Semiparametric Estimation}, pdfauthor={Juan Carlos Escanciano}, pdfkeywords={Two-step estimation; Stochastic Equicontinuity;
Semiparametric inference} }
\begin{document}

\title{Conditional Influence Functions\thanks{chernozhukov@gmail.com, wnewey@mit.edu, vsyrgk@stanford.edu

The present paper formed the basis of
an invited talk by Whitney Newey at the 2024 ESIF Conference on Economics and
AI+ML. This research was supported by NSF Grant 224247. Helpful comments were
provided by M. Kolesar and R. Singh.}}
\author{Victor Chernozhukov\\\textit{MIT}
\and Whitney K. Newey\\\textit{MIT and NBER}
\and Vasilis Syrgkanis\\\textit{Stanford University}}
\date{November 2024}
\maketitle

\begin{abstract}
There are many nonparametric objects of interest that are a function of a
conditional distribution. One important example is an average treatment effect
conditional on a subset of covariates. Many of these objects have a
conditional influence function that generalizes the classical influence
function of a functional of a (unconditional) distribution. Conditional
influence functions have important uses analogous to those of the classical
influence function. They can be used to construct Neyman orthogonal estimating
equations for conditional objects of interest that depend on high dimensional
regressions. They can be used to formulate local policy effects and describe
the effect of local misspecification on conditional objects of interest. We
derive conditional influence functions for functionals of conditional means
and other features of the conditional distribution of an outcome variable. We
show how these can be used for locally linear estimation of conditional
objects of interest. We give rate conditions for first step machine learners
to have no effect on asymptotic distributions of locally linear estimators. We
also give a general construction of Neyman orthogonal estimating equations for
conditional objects of interest.


JEL\ Classification: C13; C14; C21; D24

\end{abstract}

\section{Introduction}

There are many nonparametric, conditional objects of interest. An example is
the conditional average treatment effect given one or more key covariates.
This object quantifies how the average effect of some treatment varies with a
few covariates. An economic example is average equivalent variation for price
changes conditional on income. This object quantifies how average welfare
effects of price changes vary by income. Such conditional objects can provide
useful guidance for policy. They also give low dimensional summaries of high
dimensional unknown functions that will often have causal and/or economic interpretations.

This paper gives conditional influence functions for conditional objects of
interest and describes how they can be used. The conditional influence
function, when it exists, is a Gateaux derivative of a conditional object with
respect to the conditional CDF of a single data observation. The conditional
influence function is exactly analogous to the classic influence function of
Hampel (1974) and Huber (1981), which is a Gateaux derivative of an object of
interest with respect to the (unconditional) CDF of a single data observation.
The conditional influence function (CIF) differs in being a Gateaux derivative
with respect to a conditional CDF rather than unconditional CDF. We show in
this paper that CIF's exist for objects that are certain functionals of
conditional means and other features of conditional distributions. We also
show that CIF's have important uses.

CIF's can be used in some ways that are exactly analogous to the classical
influence function. They can be used to form Neyman orthogonal versions of
conditional moment functions that identify conditional objects of interest.
Such Neyman orthogonal conditional moment functions can be used to construct
estimators of conditional objects of interest from machine learners of
features of conditional distributions. The resulting estimators have the
property that first step machine learning has no first order effect on the
objects of interest. Consequently, bias from regularization and model
selection has little effect on objects of interest. We give a wide range of
Neyman orthogonal moment functions for objects that depend on features of the
conditional distribution of an outcome variable given regressors. We also give
a general orthogonality result that is local with respect to unknown functions
that identify the object of interest and global with respect to other unknown
debiasing functions.

The CIF also has other uses that are analogous to those of the classic
influence function. We show that it can be used to characterize local effects
of misspecification on conditional objects of interest. For example we
quantify how endogeneity of prices may affect average equivalent variation
according to income level. We also show how CIF's can be used to construct
local policy effects that are conditional on observed variables. These effects
provides ways of quantifying how effects of policy changes vary with observed
characteristics of populations, like income.

It seems CIF's cannot be used to characterize asymptotic variances or make
efficiency comparisons for objects of interest that condition on continuous
variables. It is important to note that, unlike in the unconditional case, the
CIF does not fully characterize the large-sample variances of estimators for
objects of interest that condition on continuous variables. Consequently, the
CIF does not determine efficiency rankings. For instance, series and kernel
estimators of conditional means exhibit different large-sample variances due
to variations in localization methods and approaches used to mitigate
localization biases. This distinction sets the conditional case apart from the
unconditional case and raises intriguing open questions, particularly about
the comparison of estimators in this case.\footnote{See Khan and Nekipelov
(2010) for some interesting steps in this direction.}

This paper mostly focuses on using Neyman orthogonal estimating equations for
debiased machine learning of conditional objects of interest. We estimate
conditional objects with machine learners of conditional expectations and
features of the conditional distribution of an outcome variable given
regressors. The estimation is to solve the conditional mean zero equation of
the CIF for the object of interest, plug-in first step machine learners for
unknown functions, and use nonparametric regression to estimate a conditional
mean. We focus on locally linear nonparametric regression in this paper. We
give conditions for these estimators to have standard nonparametric,
asymptotic distributions where first step machine learning can be ignored.

Estimators of conditional objects based on Neyman orthogonal moment functions
have previously been developed. Convergence rates for some conditional
objects of interest are given in Foster and Syrgkanis (2019). Linear projection
estimators of conditional objects of interest were given in Semenova and
Chernozhukov (2021). Kernel estimators with automatic debiasing are given in
Chernozhukov, Newey, and Singh (2022a). Estimation of the average treatment
effect conditional on all covariates was considered in Kennedy (2023) and
Kennedy et al. (2024). Here we give results that only require mean-square convergence rates for first step machine learners. Leudtke (2024) gave attainable rates of convergence and
inference for conditional objects of interest in reducing kernel spaces and
showed pathwise differentiability of some objects in more general spaces. Here
we consider pathwise differentiability with respect to the conditional
distribution, a different localization, and also consider many objects of
interest not covered in Leudtke (2024).

In Section 2 we define the CIF and give a wide variety of conditional objects
of interest that are linear functionals of a conditional mean. Section 3 gives
nonparametric estimators of objects of interest focusing on locally linear
regression. Section 4 describes some other uses of CIF's. Section 5 gives
CIF's for objects of interest that depend on a general feature of the
conditional distribution of an outcome given regressors. Section 6 gives a
general construction of Neyman orthogonal estimating equations for conditional parameters.

\section{The Conditional Influence Function and Examples}

Throughout the paper we will focus on the setting where data are i.i.d. with
one data observation denoted by $W.$ There are many objects of interest that
are conditional on some observed variable $V$ that is a function of $W.$ We
will use as a running example an average treatment effect conditional on some
function of covariates.

\bigskip

\textsc{Example 1} (Running): Here $W=(Y,D,Z)$ where $Y$ is an outcome,
$D\in\{0,1\}$ a treatment indicator, and $D$ is mean independent of potential
outcomes conditional on $Z$. Let $V=T(Z)$ be some function of the covariates
and $\gamma_{0}(D,Z)=E[Y|D,Z]$. Consider an object of interest%
\[
\theta_{0}(V)=E[\gamma_{0}(1,Z)-\gamma_{0}(0,Z)|V].
\]
This $\theta_{0}(V)$ is a conditional average treatment effect (CATE) given
$V$ when $0<\Pr(D=1|Z)<1$ with probability one. For example $V$ could be
family size, in which case $\theta_{0}(V)$ would quantify how the average
treatment effect depends on family size.

\bigskip

In this paper we consider objects of interest that are functions of the
conditional CDF of $W$ given $V.$ Let $F^{V}$ denote such a conditional CDF
and $F_{0}^{V}$ denote its true value, that is the conditional CDF\ of an
observed data observation $W$ given $V.$ We consider objects that take the
form%
\[
\theta_{0}(V)=\theta(F_{0}^{V}),
\]
where $\theta(F^{V})$ is a mapping from conditional CDF's of $W$ given $V$ to
the real line. The CIF will be the Gateaux derivative of $\theta(F^{V})$ with
respect to $F^{V}$ at $F^{V}=F_{0}^{V}.$

\bigskip

\textsc{Example 1} Continued: To formulate the CATE $\theta_{0}(V)$ as a
function of the CDF\ $F^{V}$ it is important that $E[Y|D,Z]$ can be viewed as
a function of $F^{V}.$ Because this formulation is of some general interest
for conditional objects other than the CATE, we give a result here:

\bigskip

\textsc{Lemma 1: }If $E[Y^{2}]<\infty$ and $V=T(X)$ for some measurable
function $T(X)$ then if $E[\{Y-\gamma(X)\}^{2}|V]$ has a unique minimizer,%
\[
E[Y|X]=\arg\min_{\gamma}E[\{Y-\gamma(X)\}^{2}|V].
\]
Proof: Note that for $\gamma_{0}(X)=E[Y|X]$ and any other function measurable
function $\gamma(X)$ with finite second moment,%
\[
E[\{Y-\gamma_{0}(X)\}^{2}|X]\leq E[\{Y-\gamma(X)\}^{2}|X].
\]
Then by iterated expectations and the above inequality,%
\begin{align*}
E[\{Y-\gamma_{0}(X)\}^{2}|V]  &  =E[E[\{Y-\gamma_{0}(X)\}^{2}|X]|V]\\
&  \leq E[E[\{Y-\gamma(X)\}^{2}|X]|V]=E[\{Y-\gamma(X)\}^{2}|V].
\end{align*}
That is,%
\begin{equation}
\gamma_{0}=\arg\min_{\gamma}E[\{Y-\gamma(X)\}^{2}|V]=\arg\min_{\gamma}%
\int\{Y-\gamma(X)\}^{2}F_{0}^{V}(dw). \label{con func}%
\end{equation}
where the minimizer is unique. $Q.E.D.$

\bigskip

By the conclusion of Lemma 1 we see that the conditional mean $E[Y|X]$ is a
function of the conditional distribution of $Y$ given $V$ equal to the
minimizing value of $E[\{Y-\gamma(X)\}^{2}|V].$ This result is well known for
$V$ equal to $X$ or constant $V$, so Lemma 1 just extends those well known
results to other functions $V$ of $X.$ This result is useful in formulating
conditional objects of interest that depend on $E[Y|X]$ as functions of the
conditional CDF $F^{V},$ as we can illustrate with Example 1.

\bigskip

\textsc{Example 1} Continued: From Lemma 1 it follows that the conditional
mean of $Y$ given $(D,Z)$ is the argmin function in equation (\ref{con func}).
Let $\gamma_{F^{V}}(D,Z)$ be this argmin function. Then we have%
\[
\theta_{0}(V)=\theta(F_{0}^{V}),\text{ }\theta(F^{V})=E_{F^{V}}[\gamma_{F^{V}%
}(1,Z)-\gamma_{F^{V}}(0,Z)|V].
\]
Thus we see that CATE is a function of the conditional CDF $F^{V}$ of $W$
given $V.$

\bigskip

We define the conditional influence function to be the Gateaux derivative of
$\theta(F^{V})$ with respect to $F^{V}$ at $F_{0}^{V}$ in a way that is
exactly analogous to the Gateaux derivative for the classic influence
function. Let $H^{V}$ be a conditional CDF that may be different than
$F_{0}^{V}$ and let
\[
F_{\tau}^{V}:=(1-\tau)F_{0}^{V}+\tau H^{V},\text{ }\tau\in\lbrack0,1],
\]
be a linear in $\tau$ deviation of the conditional CDF of $W$ given $V$ away
from its true value $F_{0}^{V}.$ This $F_{\tau}^{V}$ generalizes the one
dimensional parametric family used to define the classic influence function to
allow $F_{\tau}^{V}$ to be conditional rather than unconditional.

\bigskip

\textsc{Definition:} \textit{The conditional influence function (when it
exists) is }$\psi_{0}^{V}(W)$\textit{\ satisfying}%
\begin{equation}
\left.  \frac{d\theta(F_{\tau}^{V})}{d\tau}\right\vert _{\tau=0}=\int\psi
_{0}^{V}(w)dH^{V}(w),\text{ }E[\psi_{0}^{V}(W)|V]=0, \label{Cond Gat Der}%
\end{equation}
\textit{for all }$H^{V}\in\mathcal{H}^{V}$\textit{ for some }$\mathcal{H}%
^{V}.$

\textit{\bigskip}

This definition of a conditional influence function generalizes the definition
of the classical influence function to allow $F^{V}$ and $H^{V}$ to be
conditional CDF's rather than unconditional ones where $V$ is constant.

Before discussing the CIF and its properties further we derive the CIF in our
running example and give other examples.

\bigskip

\textsc{Example 1} continued: To derive the CIF for CATE let $E_{\tau}%
[\cdot|V]$ denote the conditional expectation given $V$ for $F_{\tau}^{V}$ and
$\gamma_{\tau}(X)=\gamma_{F_{\tau}^{V}}(X).$ Then by the chain rule of
calculus, for $\Delta_{\tau}(Z):=\gamma_{\tau}(1,Z)-\gamma_{\tau}(0,Z),$%
\[
\frac{d\theta(F_{\tau}^{V})}{d\tau}=\frac{d}{d\tau}E_{\tau}[\Delta_{\tau
}(Z)|V]=\frac{d}{d\tau}E_{\tau}[\Delta_{0}(Z)|V]+\frac{d}{d\tau}E[\Delta
_{\tau}(Z)|V].
\]
where the derivatives are evaluated at $\tau=0$ throughout. By the form of
$F_{\tau}^{V}$ it follows that%
\begin{equation}
\frac{d}{d\tau}E_{\tau}[\Delta_{0}(Z)|V]=\int[\Delta_{0}(z)-\theta
_{0}(V)]H^{V}(dw). \label{own cif}%
\end{equation}
Also for the propensity ,score $\pi_{0}(Z)=\Pr(D=1|Z)$ it follows by $V$ a
function of $Z$ and iterated expectations that
\begin{align*}
E[\gamma_{\tau}(1,Z)|V]  &  =E[E[\frac{D}{\pi_{0}(Z)}|Z]\gamma_{\tau
}(1,Z)|V]=E[\frac{D}{\pi_{0}(Z)}\gamma_{\tau}(X)|V],\\
E[\gamma_{\tau}(0,Z)|V]  &  =E[\frac{1-D}{1-\pi_{0}(Z)}\gamma_{\tau}(X)|V],
\end{align*}
so that
\[
E[\Delta_{\tau}(Z)|V]=E[\alpha_{0}(X)\gamma_{\tau}(X)|V].
\]
Then proceeding as in Newey (1994), with $E_{\tau}[\cdot|V]$ replacing
$E_{\tau}[\cdot],$ another application of the chain rule gives%
\begin{align}
\frac{d}{d\tau}E[\Delta_{\tau}(Z)|V]  &  =\frac{d}{d\tau}E[\alpha_{0}%
(X)\gamma_{\tau}(X)|V]\label{rrcif}\\
&  =\frac{d}{d\tau}E_{\tau}[\alpha_{0}(X)\gamma_{\tau}(X)|V]-\frac{d}{d\tau
}E_{\tau}[\alpha_{0}(X)\gamma_{0}(X)|V]\nonumber\\
&  =\frac{d}{d\tau}E_{\tau}[\alpha_{0}(X)Y|V]-\frac{d}{d\tau}E_{\tau}%
[\alpha_{0}(X)\gamma_{0}(X)|V]\nonumber\\
&  =\frac{d}{d\tau}E_{\tau}[\alpha_{0}(X)\{Y-\gamma_{0}(X)\}|V]\nonumber\\
&  =\int\alpha_{0}(x)[y-\gamma_{0}(x)]H^{V}(dw),\nonumber
\end{align}
where the third equality follows by iterated expectations and $V$ being a
function of $X$ and the fifth equality as in equation (\ref{own cif}).
Combining this equation with equation (\ref{own cif}) it follows that the CIF
of CATE is
\[
\psi_{0}(W)=\gamma_{0}(1,Z)-\gamma_{0}(0,Z)-\theta_{0}(V)+\alpha
_{0}(X)[Y-\gamma_{0}(X)].
\]

\bigskip

There are many other objects of interest that have a CIF. A general class of
such objects is linear functions of the conditional mean $E[Y|X]$ that are
conditionally mean square continuous. To describe these objects let $\gamma$
represent a function of $X$ that is a possible conditional mean $E[Y|X]$ and
$m(W,\gamma)$ be a linear functional of $\gamma.$ Consider any conditional
object of interest that has the form%
\[
\theta_{0}(V)=E[m(W,\gamma_{0})|V].
\]
Example 1, CATE, is a special case with $m(W,\gamma)=\gamma(1,Z)-\gamma(0,Z).$
Suppose that $m(W,\gamma)$ satisfies the following condition.

\bigskip

\textsc{Assumption 1:} \textit{(Conditional Mean Square Continuity) There
exists }$\alpha_{0}(X)$\textit{ such that }$E[\alpha_{0}(X)^{2}]<\infty
$\textit{ and}
\begin{equation}
E[m(W,\gamma)|V]=E[\alpha_{0}(X)\gamma(X)|V]\text{ for all }\gamma\text{ with
}E[\gamma(X)^{2}]<\infty. \label{conrr}%
\end{equation}

\bigskip

By a conditional version of the Riesz representation theorem existence of
$\alpha_{0}(X)$ such that equation (\ref{conrr}) is satisfied is equivalent to
$\left\vert E[m(W,\gamma)|V]\right\vert \leq D(V)E[\gamma(X)^{2}|V]$ for some
$D(V)<\infty$ with $E[\alpha_{0}(X)^{2}]<\infty.$ We refer to this $\alpha
_{0}(X)$ as a conditional Riesz representer. Under this condition the object
of interest $\theta_{0}(V)=E[m(W,\gamma_{0})|V]$ has a CIF analogous to the
CATE running example.

\bigskip

\textsc{Proposition 1:} \textit{If Assumption 1 is satisfied then the CIF of
}$\theta_{0}(V)=E[m(W,\gamma_{0})|V]$\textit{ exists and is equal to}%
\[
\psi_{0}(W)=m(W,\gamma_{0})-\theta_{0}(V)+\alpha_{0}(X)[Y-\gamma_{0}(X)].
\]

\bigskip

Proof: Let $E_{\tau}[\cdot|V]$ denote the conditional expectation given $V$
for $F_{\tau}^{V}$ and $\gamma_{\tau}(X)=\gamma_{F_{\tau}^{V}}(X)$ as in Lemma
1, so that $\theta_{F_{\tau}^{V}}(V)=E_{\tau}[m(W,\gamma_{\tau})|V].$ Then by
the chain rule it follows as in equations (\ref{own cif}) and (\ref{rrcif})
that%
\begin{align*}
\frac{d\theta(F_{\tau}^{V})}{d\tau}  &  =\frac{d}{d\tau}E_{\tau}%
[m(W,\gamma_{\tau})|V]=\frac{d}{d\tau}E_{\tau}[m(W,\gamma_{0})|V]+\frac
{d}{d\tau}E[m(W,\gamma_{\tau})|V]\\
&  =\int[m(w,\gamma_{0})-\theta_{0}(V)]H^{V}(dw)+\frac{d}{d\tau}E[\alpha
_{0}(X)\gamma_{\tau}(X)]\\
&  =\int[m(w,\gamma_{0})-\theta_{0}(V)]H^{V}(dw)+\frac{d}{d\tau}E_{\tau
}[\alpha_{0}(X)\{Y-\gamma_{0}(X)\}]\\
&  =\int\psi_{0}(w)H^{V}(dw).
\end{align*}

Also, by iterated expectations,%
\begin{align*}
E[\psi_{0}(W)|V]  &  =E[m(W,\gamma_{0})|V]-\theta_{0}(V)+E[\alpha
_{0}(X)\{Y-\gamma_{0}(X)\}|V]\\
&  =0+E[E[\alpha_{0}(X)\{Y-\gamma_{0}(X)\}|X]|V]=0.\text{ \ \ }Q.E.D.
\end{align*}

\bigskip

There are many examples of conditional objects of interest where Assumption 1
is satisfied. Here are two others in addition to the CATE.

\bigskip

\textsc{Example 2:} (CATE for continuous treatment): A second example is a
CATE for continuous treatment where $D$ is a continuous treatment variable
and
\[
\theta_{0}(V)=E[\frac{\partial}{\partial d}\gamma_{0}(D,Z)|V],\text{ }V=T(Z).
\]
As shown in Imbens and Newey (2009) the unconditional version of this object,
where $V$ is constant, is an average treatment effect when $D$ is independent
of potential outcomes conditional on $Z.$ The $\theta_{0}(V)$ given here is a
conditional version, giving the average treatment effect conditional on a
function $V$ of the covariates $Z$. It follows by integration by parts that,
when the conditional pdf $f(D|Z)$ is zero at the boundary of the support of
$D$ conditional on $Z,$ Assumption 1 is satisfied with%
\[
m(W,\gamma)=\partial\gamma(D,Z)/\partial d,\text{ }\alpha_{0}(X)=-\partial\ln
f(D|Z)/\partial d.
\]

\bigskip

\textsc{Example 3:} (Average Equivalent Variation Bound) An explicit economic
example is a bound on average equivalent variation for heterogenous demand.
Here $Y$ is the share of income spent on a commodity and $X=(P_{1},Z),$ where
$P_{1}$ is the price of the commodity and $Z$ includes income $Z_{1}$, prices
of other goods, and other observable variables affecting utility. For
$\check{P}_{1}(Z)<\bar{P}_{1}(Z)$ being lower and upper prices over which the
price of the commodity can change, $\kappa$ a bound on the income effect, and
$\omega(Z)$ some weight function. The object of interest is
\[
\theta_{0}(V)=E[\omega(Z)\int_{\check{P}_{1}}^{\bar{P}_{1}}\left(  \frac
{Z_{1}}{u}\right)  \gamma_{0}(u,Z)\exp(-\kappa\lbrack u-\check{P}%
_{1}])du|V],\text{ }V=T(Z),
\]
where $u$ is a variable of integration. When individual heterogeneity in
consumer preferences is independent of $X$ and $\kappa$ is a lower (upper)
bound on the derivative of consumption with respect to income across all
individuals, then $\theta_{0}$ is an upper (lower) bound on the weighted
average over consumers of equivalent variation for a change in the price of
the first good from $\check{P}_{1}$ to $\bar{P}_{1}$; see Hausman and Newey
(2016). Here $m(W,\gamma)=\omega(Z)\int_{\check{P}_{1}}^{\bar{P}_{1}}%
(Z_{1}/u)\gamma(u,Z)\exp(-\kappa\lbrack u-\check{P}_{1}])du$ and%
\[
\alpha_{0}(X)=f(P_{1}|Z)^{-1}\omega(Z)1(\check{P}_{1}<P_{1}<\bar{P}_{1}%
)(Z_{1}/P_{1})\exp(-\kappa\lbrack P_{1}-\check{P}_{1}]),
\]
where $f(P_{1}|Z)$ is the conditional pdf of $P_{1}$ given $Z.$

\bigskip

In each of Examples 1-3 the conditioning variable $V$ is not a one-to-one
function of $X.$ Also, the conditional Riesz representation comes from
averaging over other parts of $X.$ In Examples 1 and 2 the conditional Riesz
averages over the treatment variable $D$ and in Example 3 from averaging over
the variable $P_{1}$. If $V$ is a one-to-one function of $X$ then the
conditional Riesz representation requires that $E[m(W,\gamma)|X]=\alpha
_{0}(X)\gamma(X),$ for which $m(W,\gamma)=A(W)\gamma(X)$ for some function
$A(W)$ is the only example of which we are aware.

It is interesting to note that the CIF of Proposition 1 is identical to the
influence function of the unconditional parameter $\theta_{0}=E[m(W,\gamma
_{0})]$ except that $\theta_{0}$ is replaced by $\theta_{0}(V)$. This fact is
consistent with by a conditional Riesz representation being stronger than
(i.e. implying) an unconditional one, which follows by iterated expectations.

\section{Conditional Debiased Machine Learning via Locally Linear Regression}

The CIF can be used to construct conditionally Neyman orthogonal moment
equations for estimating conditional objects of interest, analogous to the use
of influence functions in estimating parameters. The conditional mean zero
property $E[\psi_{0}(W)|V]=0$ can be used to estimate $\theta_{0}(V)$ by
replacing $\psi_{0}(W)$ with an estimated CIF and solving for $\theta_{0}(V)$
from a nonparametric regression. For simplicity we focus in this Section on
the CIF of a linear function of a conditional mean and give a general analysis
in Section 6.

Solving $E[\psi_{0}(W)|V]=0$ in Proposition 1 for $\theta_{0}(V)$ gives%
\[
\theta_{0}(V)=E[m(W,\gamma_{0})+\alpha_{0}(X)\{Y-\gamma_{0}(X)\}|V].
\]
We can use this equation to estimate $\theta_{0}(V)$ by a nonparametric
regression of $m(W,\hat{\gamma})+\hat{\alpha}(X)\{Y-\hat{\gamma}(X)\}$ on $V$
where $\hat{\gamma}$ and $\hat{\alpha}$ are preliminary machine learners. That
is, we use $m(W,\hat{\gamma})+\hat{\alpha}(X)\{Y-\hat{\gamma}(X)\}$ as the
outcome variable in a nonparametric regression. The first part $m(W,\hat
{\gamma})$ of this outcome variable identifies $\theta_{0}(V)=E[m(W,\gamma
_{0})|V].$ The second part $\hat{\alpha}(X)\{Y-\hat{\gamma}(X)\}$ acts as a
bias correction. This second part makes the conditional expectation
$E[m(W,\gamma)+\alpha(X)\{Y-\gamma(X)\}|V]$ Neyman orthogonal with respect to
$\gamma$ and $\alpha$, thus mitigating regularization and/or model selection
bias from machine learning.

Neyman orthogonality with respect to $\gamma$ and $\alpha$ follows from the
general analysis of Section 6 and can also be shown directly. We include that
demonstration here for completeness. By using iterated expectations and
collecting terms we have%
\begin{align*}
&  E[m(W,\gamma)+\alpha(X)\{Y-\gamma(X)\}-m(W,\gamma_{0})-\alpha
_{0}(X)\{Y-\gamma_{0}(X)\}|V]\\
&  =E[m(W,\gamma-\gamma_{0})+\alpha(X)\{\gamma_{0}(X)-\gamma(X)\}|V]\\
&  =E[\alpha_{0}(X)\{\gamma(X)-\gamma_{0}(X)\}-\alpha(X)\{\gamma(X)-\gamma
_{0}(X)\}|V]\\
&  =-E[\{\alpha(X)-\alpha_{0}(X)\}\{\gamma(X)-\gamma_{0}(X)\}|V].
\end{align*}
We see that $\gamma$ and $\alpha$ being different than $\gamma_{0}$ and
$\alpha_{0}$ has zero first order effect with an explicit second order
remainder following the last equality.

To avoid over-fitting bias and keep regularity conditions for $\hat{\alpha}$
and $\hat{\gamma}$ as weak as we can cross-fitting can be used in construction
of $m(W,\hat{\gamma})+\hat{\alpha}(X)\{Y-\hat{\gamma}(X)\}.$ For cross-fitting
we partition the observation indices $i=1,...,n$ into distinct sets $I_{\ell
},$ $(\ell=1,...,L),$ and let $\hat{\gamma}_{\ell}$ and $\hat{\alpha}_{\ell}$
be computed from observations not in $I_{\ell}.$ Then for each $i\in I_{\ell}$
we take
\begin{equation}
\hat{S}_{i}=m(W_{i},\hat{\gamma}_{\ell})+\hat{\alpha}_{\ell}(X_{i})[Y_{i}%
-\hat{\gamma}_{\ell}(X_{i})] \label{debiiased}%
\end{equation}
These $\hat{S}_{i}$ observations are debiased outcomes that can be used to
estimate $\theta_{0}(V)$ from a nonparametric regression of $\hat{S}_{i}$ on
$V_{i}.$ This can be done by standard kernel or series regression or by some
machine learner if the dimension of $V$ is not small. We focus here on locally
linear regression which has good properties as a nonparametric estimator for
low dimensional $V.$ Locally linear regression is unbiased when the
$\theta_{0}(V)$ is linear in $V$, adapts well to boundary of the support of
$V$, and has improved mean square error relative to kernel regression, as
discussed in Fan (1993).

To describe a locally linear regression estimator of $\theta_{0}(V)$ based on
the debiased outcomes $\hat{S}_{i},$ let $K(u)$ be a kernel with $\int
K(u)du=1$ and for bandwidth $h>0$ and $V$ with dimension $r$ let
$K_{h}(u)=h^{-r}K(u/h).$ A locally linear estimator is%
\begin{equation}
\hat{\theta}(v)=\arg\min\limits_{\theta,\beta}\sum_{\ell=1}^{L}\sum_{i\in
I_{\ell}}(\hat{S}_{i}-\theta-(v-V_{i})^{\prime}\beta)^{2}K_{h}(v-V_{i}).
\label{llest}%
\end{equation}
This $\hat{\theta}(v)$ can be computed as the estimator of the coefficient of
$1$ in a weighted least squares regression of $\hat{S}_{i}$ on on
$(1,v-V_{i})$ with weight $K_{h}(v-V_{i})$ for the $i^{th}$ observation. Under
regularity conditions given later in this Section this $\hat{\theta}(v)$ will
have an asymptotically normal distribution that is the same as a locally
linear regression with outcome variable $S_{i}=m(W_{i},\gamma_{0})+\alpha
_{0}(X_{i})[Y_{i}-\gamma_{0}(X_{i})].$ In this sense the presence of the
estimators $\hat{\alpha}$ and $\hat{\gamma}$ can be ignored for the purposes
of choosing the bandwidth $h$ and construction of standard errors. Thus, the
bandwidth choice can be made and confidence intervals formed as in standard
algorithms for locally linear regression; see e.g. Fan and Gijbels (1996).

The bias correction term \(\alpha_{0}(X_{i})[Y_{i}-\gamma_{0}(X_{i})]\) appears in the asymptotic variance of this locally linear estimator. Moreover, the same term also appears in the variance if we use traditional (asymptotically linear) series or kernel estimators for \(\gamma\) (which is feasible when \(X\) is low-dimensional) and simply use the plug-in \(m(W_{i}, \hat{\gamma})\) in the local regression above. Thus, in the traditional low-dimensional case, the de-biasing approach and the plug-in approach are asymptotically equivalent. However, the plug-in and de-biasing approaches cease to be equivalent in the high-dimensional case because we must use regularization and/or selection to estimate \(\gamma\) well. As a result, the plug-in approach accumulates excess bias, which is removed by the bias correction term.

%presence of this term is not specific to using
%the debiased outcomes $\hat{S}_{i}.$ This term would also show up in the
%asymptotic variance of plug-in estimators based on the nonparametric
%regression of $m(W_{i},\hat{\gamma})$ that accounted for presence of 
%$\hat{\gamma}$ in the asymptotic variance. 
%Thus, the bias correction in
%$\hat{S}_{i}$ does not lead to a change in asymptotic variance but is rather a
%forward looking procedure that brings in these terms at the beginning in order
%to bias correct for possible regularization and model selection in
%$\hat{\gamma}$ and $\hat{\alpha}.$

An estimator $\hat{\alpha}(X_{i})$ of the conditional Riesz representer
$\alpha_{0}(X)$ is needed for construction of $\hat{S}_{i}.$ As noted in
Section 3, the conditional Riesz representation of Assumption 1 implies the
unconditional Riesz representation%
\[
E[m(W,\gamma)]=E[\alpha_{0}(X)\gamma(X)]\text{ for all }\gamma\text{ with
}E[\gamma(X)^{2}]<\infty.
\]
As in Chernozhukov et al. (2024), this equation is the first-order conditions
for the extremum problem%
\[
\alpha_{0}(X)=\arg\min_{\alpha}\{E[-2m(W,\gamma)+\alpha(X)^{2}]\}.
\]
Thus $\alpha_{0}(X)$ can be estimated by minimizing a sample version of this
objective function over some approximating set $\mathcal{A}_{n}$, with
possibly a penalty $\hat{P}_{n}(\alpha)$ included in the objective, as in%
\[
\hat{\alpha}_{\ell}=\arg_{\alpha\in\mathcal{A}_{n}}\min\{\sum_{i\notin
I_{\ell}}[-2m(W_{i},\alpha)+\alpha(X_{i})^{2}]+\hat{P}_{n}(\alpha)\}.
\]
Mean square convergence rates for this $\hat{\alpha}_{\ell}$ are given for
Lasso in Chernozhukov, Newey, and Singh (2022b) and for other estimators in
Chernozhukov et al. (2024). These are automatic methods for constructing
$\hat{\alpha},$ requiring only $m(W,\alpha)$ and not relying on an explicit
formula for $\alpha_{0}.$

Under regularity condition the presence of $\hat{\gamma}$ and $\hat{\alpha}$
will not affect the nonparametric asymptotic distribution of $\hat{\theta}(v)$
because the nonparametric regression is being done on the debiased outcome
$\hat{S}_{i}$. These conditions will allow standard nonparametric inference
methods for local linear regression to be applied using $\hat{\theta}(v)$. In
particular, with an under-smoothing choice for $h$ it will be the case that $\sqrt{nh^{r}%
}[\hat{\theta}(v)-\theta_{0}(v)]$ is asymptotically normal with the same
asymptotic variance as for the true debiased outcome $S_{i}=m(W_{i},\gamma
_{0})+\alpha_{0}(X_{i})\{Y_{i}-\gamma_{0}(X_{i})\}$. Here $\sqrt{nh^{r}}$ is
the standard normalization for asymptotic normality of a locally linear regression.

To specify conditions for $\hat{\theta}(v)$ to have these useful large sample
properties we will give conditions sufficient for $\sqrt{nh^{r}}[\hat{\theta
}(v)-\tilde{\theta}(v)]\overset{p}{\longrightarrow}0,$ where $\tilde{\theta
}(v)$ is the random variable obtained by replacing $\hat{S}_{i}$ by $S_{i}$ in
equation (\ref{llest}). We impose mild conditions on the kernel $K(u).$

\bigskip

\textsc{Assumption 2:} $\int K(u)du=1$ \textit{and }$K(u)$\textit{ and
}$K(u)u^{2}$\textit{ are bounded.}

\bigskip

We also assume that $m(W,\gamma)$ is a mean square continuous function of
$\gamma$ and that the Riesz representer and $Var(Y|X)$ are bounded.

\bigskip

\textsc{Assumption 3:} $E[m(W,\gamma)^{2}]\leq CE[\gamma(X)^{2}]$ \textit{for
some }$C>0$\textit{ and} $\alpha_{0}(X)$\textit{ and }$Var(Y|X)$\textit{ are
bounded}.

\bigskip

For any function $h(X)$ let $\left\Vert h\right\Vert =\sqrt{E[h(X)^{2}]}$
denote the mean-square norm. Our main regularity conditions are mean square
rates of convergence for $\hat{\gamma}$ and $\hat{\alpha}$ relative to the
bandwidth $h$ of the kernel estimator.

\bigskip

\textsc{Assumption 4: }i)$\sqrt{nh^{r}}\longrightarrow\infty,$ $\left\Vert
\hat{\gamma}-\gamma_{0}\right\Vert =o_{p}(h^{r/2}),$ $\left\Vert \hat{\alpha
}-\alpha_{0}\right\Vert =o_{p}(h^{r/2}),$ ii) $\sqrt{n}\left\Vert \hat{\gamma
}-\gamma_{0}\right\Vert \left\Vert \hat{\alpha}-\alpha_{0}\right\Vert
=o_{p}(h^{r/2}).$

\bigskip

These conditions impose convergence rates for the mean square error of
$\hat{\gamma}$ and $\hat{\alpha}$ as well as a combined rate for the two functions.

\bigskip

\textsc{Theorem 2:} \textit{If Assumptions 1 - 4 are satisfied then}%
\[
\sqrt{nh^{r}}\{\hat{\theta}(v)-\tilde{\theta}(v)\}\overset{p}{\longrightarrow
}0.
\]

\bigskip

This result is proved in the Appendix. To explain conditions that are implicit
in Assumption 4 we give some necessary conditions for that assumption when the
mean-square convergence rates of $\hat{\gamma}$ and $\hat{\alpha}$ are
optimal. If $X$ had dimension $d,$ $X$ has compact support, and $\gamma
_{0}(X)$ has $s$ derivatives then from Stone (1980) the best attainable rate
for $\left\Vert \hat{\gamma}-\gamma_{0}\right\Vert $ is $n^{-s/(d+2s)}.$ For
simplicity suppose that both $\left\Vert \hat{\gamma}-\gamma_{0}\right\Vert
=O_{p}(n^{-s/(d+2s)})$ and $\left\Vert \hat{\alpha}-\alpha_{0}\right\Vert
=O_{p}(n^{-s/(d+2s)}),$ corresponding to $\alpha_{0}$ also having $s$
derivatives and having the same optimal rate as $\hat{\gamma}.$ Let $V$ have
dimension $r,$ and suppose that $\theta_{0}(V)$ also has bounded second
derivatives. Then a necessary condition for Assumption 3 is%
\begin{equation}
\frac{r+2}{r+4}<\frac{2s}{d+2s}. \label{rate bound}%
\end{equation}
For any positive integers $d$ and $r$ $\leq d$ this condition will be
satisfied for $s$ large enough. For example, for $r=1,$ where $V$ is a scalar,
this condition can be shown to be $s>3d/4.$ We note that this condition is
stronger than the corresponding condition $s>d/2$ for both of $\hat{\gamma}$
and $\hat{\alpha}$ to converge faster in mean square than $n^{-1/4},$ as is to
be expected from Assumption 3, which requires that $\left\Vert \hat{\gamma
}-\gamma_{0}\right\Vert \left\Vert \hat{\alpha}-\alpha_{0}\right\Vert $ goes
to zero faster than $1/\sqrt{n}.$

The conditions of Theorem 2 only require $L_2$ convergence, i.e. convergence in mean-square, of $\hat\gamma$ and $\hat\alpha$ at specified rates. 
Hypotheses of $L_2$ convergence makes this result widely applicable to machine learning, where some learners are only known to have $L_2$ rates. 
If instead $L_4$ rates were imposed, e.g. as in Foster and Srygkannis (2023), then equation (3.3) could be improved.
Further improvements may also possible, e.g. as in Kennedy (2023), for $\hat\gamma$ and $\hat\alpha$ with special structure.
We have chosen to focus here on mean-square rate conditions because they apply most widely to first step machine learners.

\section{Other Uses of the Conditional Influence Function}

Equation (\ref{Cond Gat Der}) motivates the use of the influence function for
economic applications. The Gateaux derivative $d\theta(F_{\tau}^{V})/d\tau$ is
the local effect of changing the distribution $F$ on the object $\theta
(F^{V}).$ If $\theta\left(  F^{V}\right)  $ is an economic object of interest,
such as a feature of the distribution of outcome variables, then
$d\theta(F_{\tau}^{V})/d\tau$ can be thought of as a local policy effect of
changing the distribution of the data. Equation (\ref{Cond Gat Der}) then can
be used to obtain the local policy effect from the conditional influence
function, generalizing Firpo, Fortin, and Lemeiux (2009) to conditional policy effects.

When $\theta(F^{V})$ is the probability limit of an estimator $\hat{\theta
}(V)$ we can think of $d\theta(F_{\tau}^{V})/d\tau$ as the local sensitivity
of that estimator to changes in $F^{V}$, which gives local effects of
misspecification. This use of the CIF generalizes the sensitivity analysis of
Andrews, Gentzkow, and Shapiro (2017) to apply to nonparametric, conditional
objects of interest. Quantifying local sensitivity of an estimator of a
conditional object to misspecification is another potentially important use of
the CIF. For example, it would be possible to check sensitivity to
misspecification of how economic effects vary by income by generalizing the sensitivity analysis of Ichimura and Newey (2021) to the CIF. The focus of this
paper is on debiased machine learning of conditional objects of interest and
so we defer discussion of using the CIF for local policy or misspecification analysis.

\section{Functions of Features of the Conditional Distribution}

In this Section we give estimators of objects of interest that depend on
features of the conditional distribution of $Y$ given $X$ beyond the
conditional mean. We give locally linear regression estimators based on
debiased outcomes and asymptotic theory sufficient for first step estimators
to have no effect on the nonparametric limiting distribution of the estimator.
The results illustrate how the CIF can be used to estimate objects that depend
on features of the conditional distribution of $Y$ given $X$ beyond the
conditional mean. At the time of the writing of this paper such conditional
objects had not previously been considered in the literature.

The features of the conditional distribution we consider are those $\gamma
_{0}(X)$ such that%
\begin{equation}
\gamma_{0}(X)=\arg\min_{\gamma}E[Q(W,\gamma(X))|X], \label{argmin}%
\end{equation}
where $Q(W,a)$ is a convex function of the scalar $a.$ The objects of interest
are any $\theta_{0}(V)$ with $\theta_{0}(V)=E[m(W,\gamma_{0})|V],$ where
$m(W,\gamma)$ is a function of a data observation and a linear function of
$\gamma.$ These objects go beyond those that depend on the conditional mean by
allowing $\gamma_{0}(X)$ to be some other feature of the conditional
distribution of $Y$ given $X$.

An example of $\theta_{0}(V)$ is a conditional average derivative of a
conditional quantile.

\bigskip

\textsc{Example 4: }For $0<\nu<1$ and $q_{v}(u)=[1(u<0)(1-\tau)+\nu
1(u>0)]\left\vert u\right\vert $ and $Q(W,\gamma)=q_{\nu}(Y-\gamma(X)).$ The
$\nu^{th}$ conditional quantile of the conditional distribution of $Y$ given
$X$ is the minimizer in equation (\ref{argmin}). Suppose that $X=(D,Z)$ where
$D$ is a continuous variable of interest and $Z$ are covariates. A conditional
object of interest is
\[
\theta_{0}(V)=E[\frac{\partial\gamma_{0}(D,Z)}{\partial d}|V],\text{ }V=T(Z).
\]
This $\theta_{0}(V)$ is an average derivative of a conditional quantile of $Y$
given $X$ conditional on a function $V$ of covariates. It is conditional on
$V$ version of the object considered by Chaudhuri, Doksum, and Samarov (1997).

\bigskip

For the CIF of $\theta_{0}(V)$ to exist the function $\gamma_{0}$ should be a
function of the conditional distribution of $W$ given $V.$ This property
follows from equation (\ref{argmin}) similarly to Lemma 1 because $\gamma
_{0}=\arg\min_{\gamma}E[Q(W,\gamma)|V]$ by iterated expectations. For brevity
we omit a formal statement of this result and just note that $\gamma
_{0}=\gamma(F_{0}^{V})$ for $\gamma(F_{V})=\arg\min_{\gamma}E_{F^{V}%
}[Q(W,\gamma)|V].$ Because $\gamma_{0}$ can be viewed as a function of $F^{V}$
the object of interest can also be viewed in this way with $\theta
_{0}(V)=\theta(F_{0}^{V})$ for%
\begin{equation}
\theta(F^{V})=E_{F^{V}}[m(W,\gamma(F^{V}))|V],\text{ }\gamma(F^{V})=\arg
\min_{\gamma(X)}E_{F^{V}}[Q(W,\gamma)|V]. \label{func}%
\end{equation}
This $\theta(F^{V})$ will have a CIF under some regularity conditions. The CIF
will depend on the derivative $\rho(W,\gamma_{0}(X))$ of $Q(W,\gamma
_{0}(X)+a)$ with respect to the constant $a$ at $a=0$ which we assume exists
with probability one. The first order condition for $\gamma_{0}(X)$ is%
\[
E[\rho(W,\gamma_{0}(X))|X]=0.
\]

We use regularity conditions in terms of $F_{\tau}^{V}=(1-\tau)F_{0}^{V}+\tau
H^{V}$ similarly to Ichimura and Newey (2022). Let $\gamma_{\tau}%
=\gamma(F_{\tau}^{V})$ and $E_{\tau}[\cdot|V]$ denote the conditional
expectation with respect to $F_{\tau}^{V}.$

\bigskip

\textsc{Assumption 5: }\textit{i) There is }$v_{m}(X)$\textit{\ such that
}$E[m(W,b)|V]=E[v_{m}(X)b(X)|V]$ for all bounded $b(X)$\textit{ and }%
$E[v_{m}(X)^{2}]<\infty;$\textit{ ii) there is }$v_{\rho}(X)<0$\textit{ that
is bounded and bounded away from zero such that }$\partial E[b(X)\rho
(W,\gamma_{\tau})|V]/\partial\tau=\partial E[b(X)v_{\rho}(X)\gamma_{\tau
}(X)|V]/\partial\tau$\textit{ for every bounded }$b(X)$\textit{.}

\bigskip

In ii) $v_{\rho}(X)$ will be the derivative of $E[\rho(W,\gamma_{0}(X)+a)|X]$
with respect to the scalar $a$ evaluated at $a=0.$ This condition allows for
$\rho(W,\gamma)$ to not be continuous as long as $E[\rho(W,\gamma
_{0}(X)+a)|X]$ is differentiable in $a$. Also $v_{\rho}(X)<0$ is a sign
normalization (that holds when $\rho(W,\gamma(W))=Y-\gamma(X)$) while
$v_{\rho}(X)$ being bounded and bounded away from zero is important for the results.

\bigskip

\textsc{Proposition 3:} \textit{If Assumption 5 is satisfied then for }%
$\alpha_{0}(X)=-v_{m}(X)/v_{\rho}(X),$ $\theta(F^{V})=E_{F^{V}}[m(W,\gamma
(F_{V}))]$ has CIF%
\[
\psi_{0}(W)=m(W,\gamma_{0})-\theta_{0}(V)+\alpha_{0}(X)\rho(W,\gamma_{0}).
\]

\bigskip

Proof: For $F_{\tau}^{V}=(1-\tau)F_{0}^{V}+\tau H^{V},$ $0<\tau<1,$ let
$\gamma_{\tau}=\gamma(F_{\tau}^{V})$ and $E_{\tau}[\cdot|V]=E_{F_{\tau}^{V}%
}[\cdot|V]$ as in Section 2. By the chain rule of calculus,%
\begin{align}
\frac{d}{d\tau}\theta(F_{\tau}^{V})  &  =\frac{d}{d\tau}E_{\tau}%
[m(W,\gamma_{\tau})|V]\label{gen der}\\
&  =\frac{d}{d\tau}E_{\tau}[m(W,\gamma_{0})|V]+\frac{d}{d\tau}E[m(W,\gamma
_{\tau})|V]\nonumber\\
&  =\int[m(w,\gamma_{0})-\theta_{0}]H^{V}(dw)+\frac{d}{d\tau}E[m(W,\gamma
_{\tau})|V].\nonumber
\end{align}
Also, by Assumption 5 i) the first order conditions for $\gamma_{\tau}$ for
every $\tau$ give%
\[
0\equiv E_{\tau}[\alpha_{0}(X)\rho(W,\gamma_{\tau})|V]
\]
identically in $\tau$. Differentiating this identity with respect to $\tau$
and applying the chain rule gives%
\begin{align*}
0  &  =\frac{d}{d\tau}E_{\tau}[\alpha_{0}(X)\rho(W,\gamma_{0})|V]+\frac
{d}{d\tau}E[\alpha_{0}(X)\rho(W,\gamma_{\tau})|V]\\
&  =\int\alpha_{0}(x)\rho(w,\gamma_{0})H^{V}(dw)+\frac{d}{d\tau}E[\alpha
_{0}(X)\rho(W,\gamma_{\tau})|V].
\end{align*}
Solving gives $-dE[\alpha_{0}(X)\rho(W,\gamma_{\tau})|V]/d\tau=\int\alpha
_{0}(x)\rho(w,\gamma_{0})H^{V}(dw).$ Then%

\begin{align*}
\frac{d}{d\tau}E[m(W,\gamma_{\tau})|V]  &  =\frac{d}{d\tau}E[v_{m}%
(X)\gamma_{\tau}(X)|V]=\frac{d}{d\tau}E[-\alpha_{0}(X)v_{\rho}(X)\gamma_{\tau
}(X)|V]\\
&  =-\frac{d}{d\tau}E[\alpha_{0}(X)\rho(W,\gamma_{\tau})|V]=\int\alpha
_{0}(x)\rho(w,\gamma_{0})H^{V}(dw),
\end{align*}
where the first equality follows by Assumption 5 i), the second by the
definition of $\alpha_{0}(X)$, the third by Assumption 5 ii), and the fourth
by the previous equation. The conclusion follows by substituting in equation
(\ref{gen der}). $Q.E.D.$

\bigskip

This result extends Theorem 1 of Ichimura and Newey (2022) to conditional
influence functions of conditional objects that depend on $\gamma_{0}(X)$ from
equation (\ref{argmin}).

A nonparametric estimator of $\theta_{0}(V)$ can be constructed using the CIF
similarly to Section 3. Solving $E[\psi_{0}(W)|V]=0$ in Proposition 3 for
$\theta_{0}(V)$ gives%
\[
\theta_{0}(V)=E[m(W,\gamma_{0})+\alpha_{0}(X)\rho(W,\gamma_{0})|V].
\]
We can use this equation to estimate $\theta_{0}(V)$ by a nonparametric
regression of $m(W,\hat{\gamma})+\hat{\alpha}(X)\rho(W,\hat{\gamma})$ on $V$
where $\hat{\gamma}$ and $\hat{\alpha}$ are preliminary machine learners. As
in Section 3 the second part $\hat{\alpha}(X)\rho(W,\hat{\gamma})$ of this
outcome variable provides a bias correction for regularization and/or model
selection bias from machine learning. Conditional Neyman orthogonality holds
by $E[\alpha(X)\rho(W,\gamma_{0})|V]=E[\alpha(X)E[\rho(W,\gamma_{0})|X]|V]=0$
for any $\alpha(X)$, which implies%
\begin{align*}
&  E[m(W,\gamma)+\alpha(X)\rho(W,\gamma)-m(W,\gamma_{0})-\alpha_{0}%
(X)\rho(W,\gamma_{0})|V]\\
&  =E[m(W,\gamma-\gamma_{0})+\alpha_{0}(X)\{\rho(W,\gamma)-\rho(W,\gamma
_{0})\}|V]\\
&  +E[\{\alpha(X)-\alpha_{0}(X)\}\{\rho(W,\gamma)-\rho(W,\gamma_{0})\}|V].
\end{align*}
The second term in this expression is an explicit second order remainder and
the first term is also second order because $\alpha_{0}(X)=-v_{m}(X)/v_{\rho
}(X)$ implies that $\alpha_{0}(X)\{\rho(W,\gamma)-\rho(W,\gamma_{0})\}$
cancels out $m(W,\gamma-\gamma_{0})$ to first order.

A cross-fit version of the debiased outcome variable
\begin{equation}
\hat{S}_{i}=m(W_{i},\hat{\gamma}_{\ell})+\hat{\alpha}_{\ell}(X_{i})\rho
(W_{i},\hat{\gamma}_{\ell}),
\end{equation}
can be used to estimated $\theta_{0}(V)$ from a nonparametric regression of
$\hat{S}_{i}$ on $V_{i}.$ A locally linear estimator can be obtained as in
equation (\ref{llest}).

An estimator $\hat{\alpha}_{\ell}(X_{i})$ is needed to construct the debiased
outcome $\hat{S}_{i}$. We follow Chernozhukov et al. (2024) in this
construction. We suppose that there is $v_{\rho}(W)$ such that $E[v_{\rho
}(W)|X]$ is approximately $v_{\rho}(X)$ and let $\hat{v}_{\rho}(W)$ be an
estimator of $v_{\rho}(W).$ Then a cross-fit estimator of $\alpha_{0}(X)$ can
be constructed as
\begin{equation}
\hat{\alpha}_{\ell}=\arg_{\alpha\in\mathcal{A}_{n}}\min\sum_{i\notin I_{\ell}%
}\{-2m(W_{i},\alpha)-\hat{v}_{\rho}(W_{i})\alpha(X_{i})^{2}\}+\hat{P}%
_{n}(\alpha), \label{alpha obj}%
\end{equation}
where $\hat{v}_{\rho}(w)$ is an estimator of $v_{\rho}(w)$. We use $\hat
{v}_{\rho}(w)$ as a function of $w$ here to avoid the need to estimate
$E[v_{\rho}(W)|X]$ in constructing the objective function.

Here we need additional conditions to those previously given to allow for
$\rho(W,\gamma(X))$ that is not $Y-\gamma(X).$ Let $\hat{\Upsilon}_{\ell
}(w)=\{\hat{\alpha}_{\ell}(x)-\alpha_{0}\}\{\rho(W,\hat{\gamma}_{\ell}%
)-\rho(W,\gamma_{0})\}.$

\bigskip

\textsc{Assumption 6: }\textit{For each }$\ell=1,...,L$\textit{, either i)}%
\[
\sqrt{n}\int\hat{\Upsilon}_{\ell}(w)F_{0}(dw)=o_{p}(h^{r/2}),\text{ }\int%
\hat{\Upsilon}_{\ell}(w)^{2}F_{0}(dw)=o_{p}(h^{r});
\]
\textit{\ or ii) }$\sqrt{n}\left\Vert \hat{\alpha}_{\ell}-\alpha
_{0}\right\Vert \left\Vert \rho(W,\hat{\gamma}_{\ell})-\rho(W,\gamma
_{0})\right\Vert =o_{p}(h^{r/2}).$

\bigskip

Here the choice between conditions i) and ii) is allowed so that Assumption 6 can apply to first step conditional quantile estimation, where i) is satisfied but not ii).

\bigskip

\textsc{Assumption 7:} \textit{There is bounded }$v_{m}(X)$ with
$E[m_{n}(W,\gamma)|V]=E[v_{m}(X)\gamma(X)|V],$\textit{ for all }%
$E[\gamma(X)^{2}]<\infty,$\textit{ }$E[\rho(W,\gamma_{0})|X]=0,$\textit{ and
there is }$R(X,\gamma,\gamma_{0})$\textit{ such that}
\begin{align*}
E[\rho(W,\gamma)|X]  &  =v_{\rho}(X)\{\gamma(X)-\gamma_{0}(X)\}+R(X,\gamma
,\gamma_{0}),\\
E[|R(X,\gamma,\gamma_{0})|]  &  \leq C\left\Vert \gamma-\gamma_{0}\right\Vert
^{2}.
\end{align*}
\textit{where }$v_{\rho}(X)>0$\textit{ is bounded and bounded away from zero.
Also either }$R(X,\gamma,\gamma_{0})=0$\textit{ or }$\sqrt{n}\left\Vert
\hat{\gamma}-\gamma_{0}\right\Vert ^{2}=o_{p}(h^{r/2})$\textit{.}

\bigskip

With these additional conditions in place we can show $\hat{\theta}(v)$ is
asymptotically equivalent to $\tilde{\theta}(v)$ obtained by using
$S_{i}=m(W_{i},\gamma_{0})+\alpha_{0}(X_{i})\rho(W_{i},\gamma_{0})$ in place
of $\hat{S}_{i}$ in equation (\ref{llest}).

\bigskip

\textsc{Theorem 4: }\textit{If Assumptions 2, 3, 4 i), 6, and 7 are satisfied
for }$\alpha_{0}(X)=-v_{m}(X)/v_{\rho}(X)$ \textit{then}
\[
\sqrt{nh^{r}}\{\hat{\theta}(v)-\tilde{\theta}(v)\}\overset{p}{\longrightarrow
}0.
\]

\section{Neyman Orthogonal Conditional Moment Functions}

A quite general way to formulate conditional objects of interest is as the
solution to a conditional moment restriction. Let $g(W,\gamma,\theta)$ be a
function of data observation $W,$ an unknown function $\gamma,$ and a possible
possible conditional object of interest $\theta(V)$. We assume that
$\theta_{0}(V)$ solves%
\begin{equation}
0=E[g(W,\gamma_{0},\theta)|V]. \label{cond gmm}%
\end{equation}
All of the conditional objects of interest we have given thus far solve such
an equation for $g(W,\gamma,\theta)=m(W,\gamma)-\theta.$ The formula in
equation (\ref{cond gmm}) allows for $\theta_{0}(V)$ to solve a conditional
moment restriction.

Neyman orthogonal moment functions can be constructed by adding to
$g(W,\gamma,\theta)$ the CIF $\phi(W,\gamma,\alpha,\theta)$ of $E[g(W,\gamma
(F^{V}),\theta)|V],$ where $\gamma(F^{V})$ is the first step, nonparametric object that helps identify $\theta_{0}(V)$ and $\alpha$ is a first step function
additional to $\gamma.$ This CIF satisfies, for $F_{\tau}^{V}=(1-\tau
)F_{0}^{V}+\tau H^{V},$%
\begin{equation}
\frac{d}{d\tau}E[g(W,\gamma(F_{\tau}^{V}),\theta)|V]=\int\phi(w,\gamma
_{0},\alpha_{0},\theta)dH^{V}(dw),\text{ }E[\phi(W,\gamma_{0},\alpha
_{0},\theta)|V]=0, \label{corr term}%
\end{equation}
for all $H^{V}\in\mathcal{H}^{V}$.
A conditionally Neyman orthogonal estimating function is%
\[
\psi(W,\gamma,\alpha,\theta)=g(W,\gamma,\theta)+\phi(W,\gamma,\alpha,\theta).
\]
In Section 3 $\phi(W,\gamma,\alpha,\theta)=\alpha(X)[Y-\gamma(X)]$ and in
Section 5 $\phi(W,\gamma,\alpha,\theta)=\alpha(X)\rho(W,\gamma).$

The function $\psi(W,\gamma,\alpha,\theta)$ has two orthogonality properties
that we derive and explain. These properties depend on assuming that the true
value $\alpha_{0}$ of $\alpha$ can also be viewed as a function of the
conditional distribution $F_{V}$ of $W$ given $V,$ i.e. that the following
condition is satisfied.

\bigskip

\textsc{Assumption 8:} \textit{There is }$\alpha(F_{V})$\textit{ such that if
}$F_{0}^{V}=F^{V}$\textit{ then }$\alpha_{0}=\alpha(F^{V}).$

\bigskip

This condition is satisfied in each of the examples we have given. For the
conditional average treatment effect $\alpha_{0}(X)$ depends on the propensity
score $\Pr(D=1|Z)$ which is a function of $F_{V},$ by Lemma 1 and by $V$ being a function of the covariates $Z$.
Assumption 8 can also be shown to hold in the other examples.

\bigskip

\textsc{Proposition 5:} \textit{If Assumption 8 is satisfied then for all
}$F_{\tau}^{V}$\textit{ such that }$\alpha(F_{\tau}^{V})=\alpha_{0}$\textit{
for all }$\tau$\textit{ in a neighborhood of zero,}%
\[
\frac{d}{d\tau}E[\psi(W,\gamma(F_{\tau}^{V}),\alpha_{0},\theta)|V]=0.
\]
\textit{Also for all }$\alpha$\textit{ such that there exists }$\alpha(F^{V}%
)$\textit{ with }$\alpha=\alpha(F^{V})$\textit{ and }$\gamma((1-\tau
)F^{V}+\tau F_{0}^{V})=\gamma_{0}$\textit{ for all }$\tau$\textit{ small
enough, }%
\[
E[\phi(W,\gamma_{0},\alpha,\theta)]=0.
\]

\bigskip

Proof: Let $F_{\tau}^{V}$ be such that $\alpha(F_{\tau})=\alpha_{0}$ for all
$\tau$ near $0$. Then by the conditional mean zero property of $\phi
(W,\gamma(F_{\tau}^{V}),\alpha_{0},\theta)$ in equation (\ref{corr term}),
\[
E_{\tau}[\phi(W,\gamma(F_{\tau}^{V}),\alpha_{0},\theta)|V]\equiv0.
\]
Differentiating this identity with respect to $\tau,$ applying the chain rule,
and evaluating at $\tau=0$ gives%
\begin{align*}
0  &  =\int\phi(W,\gamma_{0},\alpha_{0},\theta)H^{V}(dw)+\frac{\partial
}{\partial\tau}E[\phi(W,\gamma(F_{\tau}^{V}),\alpha_{0},\theta)|V]\\
&  =\frac{d}{d\tau}E[g(W,\gamma(F_{\tau}^{V}),\theta)|V]+\frac{\partial
}{\partial\tau}E[\phi(W,\gamma(F_{\tau}^{V}),\alpha_{0},\theta)|V]\\
&  =\frac{d}{d\tau}E[\psi(W,\gamma(F_{\tau}^{V}),\alpha_{0},\theta)|V],
\end{align*}
where the second equality holds by equation (\ref{corr term}), giving the
first conclusion.

For the second conclusion note that for $\bar{F}_{\tau}^{V}=(1-\tau)F^{V}+\tau
F_{0}^{V}$, by $\gamma(F^{V})=\gamma_{0},$
\begin{align*}
E[\phi(W,\gamma_{0},\alpha,\theta)|V]  &  =\int\phi(w,\gamma_{0},\alpha
(F^{V}),\theta)F_{0}^{v}(dw)\\
&  =\frac{\partial}{\partial\tau}E_{F^{V}}[g(W,\gamma(\bar{F}_{\tau}%
^{V}),\theta)|V]\\
&  =\frac{\partial}{\partial\tau}E_{F^{V}}[g(W,\gamma_{0},\theta)|V]=0.
\end{align*}
where the second equality follows by equation (\ref{corr term}) and the third
equality by by $\gamma(\bar{F}_{\tau}^{V})=\gamma_{0}$ for all $\tau$ near
enough to zero. $Q.E.D.$

\bigskip

This result requires that $\gamma\left(  F^{V}\right)  $ depend on different
features of $F_{V}$ than $\alpha(F^{V}).$ The first conclusion comes from
$\gamma(F^{V})$ varying while $\alpha(F^{V})=\alpha_{0}$ and the second
conclusion comes from $\alpha(F^{V})$ varying while $\gamma(F^{V})=\gamma
_{0}.$ In each of the examples such variation is possible because because
$\gamma(F^{V})$ is some feature of the conditional distribution of an outcome
variable $Y$ given regressors $X$ and $\alpha(F^{V})$ is a feature of the
conditional distribution of some of the regressors $X$ conditional on $V.$ We
conjecture that this property is satisfied in many other settings also.

The second conclusion is a strong form of Neyman orthogonality of
$\psi(W,\gamma_{0},\alpha,\theta_{0})$ in $\alpha,$ implying that
\[
E[\psi(W,\gamma_{0},\alpha,\theta_{0})|V]=E[g(W,\gamma_{0},\theta
_{0})|V]+E[\phi(W,\gamma_{0},\alpha,\theta_{0})|V]=0.
\]
Here the Neyman orthogonal moment function has zero conditional expectation
for any value of $\alpha$ and not just for $\alpha_{0}.$ Note that this
property holds in the examples we have given because $E[\rho(W,\gamma
_{0})|X]=0,$ with $\rho(W,\gamma)=Y-\gamma(X)$ in Section 2. We conjecture
that this robustness property holds in many other settings also.

The first conclusion is Neyman orthogonality of $\psi(W,\gamma,\alpha
_{0},\theta_{0})$ as $\gamma$ varies away from $\gamma_{0}$ along each path of
the form $\gamma(F_{\tau}^{V}).$ Under additional regularity conditions this
pathwise orthogonality will imply orthogonality for other kinds of variation
of $\gamma$, such as Gateaux differentiability where orthogonality is%
\[
\left.  \frac{\partial}{\partial\delta}E[\psi(W,\gamma_{0}+\delta
(\gamma-\gamma_{0}),\alpha_{0},\theta_{0})|V]\right\vert _{\delta=0}=0,
\]
for all $\gamma\in\Gamma,$ where $\Gamma$ is a set of $\gamma$. For example,
one could impose conditional versions of the hypotheses of Theorem 3 of
Chernozhukov et al. (2022c) to obtain this orthogonality. 
We focus here on the pathwise orthogonality of Proposition 5 and leave the establishment of other forms of conditional orthogonality to particular conditional parameters of interest.

When $E[\psi(W,\gamma,\alpha_{0},\theta_{0})|V]$ is linear affine in $\gamma$
the first conclusion of Proposition 5 generally implies a stronger
orthogonality property that $E[\psi(W,\gamma,\alpha_{0},\theta_{0})|V]$ for
all $\gamma.$ Together with the second conclusion of Proposition 5, we see
that the conditional Neyman orthogonal moment function $\psi(W,\gamma
,\alpha,\theta)$ will be conditionally double robust, where $E[\psi
(W,\gamma,\alpha_{0},\theta_{0})|V]=E[\psi(W,\gamma_{0},\alpha,\theta
_{0})|V]=0,$when $E[\psi(W,\gamma,\alpha_{0},\theta_{0})|V]$ is linear affine
in $\gamma$.

Proposition 5 generalizes results of Section 4 of Chernozhukov et al. (2022c) to the CIF.

\pagebreak

\section{APPENDIX}

In this Appendix we will prove Theorems 4 and 6. We first give a Lemma that is the basis of Theorems 4 and 6. Let $m_{n}(W,\gamma)$ be a linear functional of
$\gamma$ that will be the product of the kernel $K_{h}(v-V)$ and the
$m(W,\gamma)$ in the body of the paper. Also, let $\alpha_{n}(X)$ and
$\hat{\alpha}_{n}(X)$ be the product of the kernel and $\alpha_{0}(X)$ and
$\hat{\alpha}(X)$ respectively. Let $\psi_{n}(w,\gamma,\alpha)=m_{n}%
(w,\gamma)+\alpha(x)\rho(w,\gamma)$ and $\sigma_{n}$ be an increasing sequence
that will be $h^{-r/2}$ for the proofs of Theorems 4 and 6. In Lemma A we give
conditions for%
\begin{equation}
\frac{1}{n}\sum_{\ell=1}^{L}\sum_{i\in I_{\ell}}\{\psi_{n}(W_{i},\hat{\gamma
}_{\ell},\hat{\alpha}_{n\ell})-\psi_{n}(W_{i},\gamma_{0},\alpha_{n}%
)\}=o_{p}(\sigma_{n}/\sqrt{n}), \label{Key}%
\end{equation}
which will give the conclusions of Theorems 4 and 6.

The first condition imposes rates of convergence for $\hat{\gamma}_{\ell}$ and
$\hat{\alpha}_{n\ell}.$

\bigskip

\textsc{Assumption A1: }$\alpha_{n}(X)$ \textit{and }$Var(\rho(W,\gamma
_{0})|X)$\textit{ are bounded and for each }$\ell=1,...,L$\textit{, i) }$\int
m_{n}(w,\hat{\gamma}_{\ell}-\gamma_{0})^{2}F_{0}(dw)=o_{p}(\sigma_{n}^{2}%
),$\textit{\ ii) }$\int\alpha_{n}(x)^{2}\{\rho(w,\hat{\gamma}_{\ell}%
)-\rho(w,\gamma_{0})\}^{2}F_{0}(dw)=o_{p}(\sigma_{n}^{2}),$\textit{\ and iii)
}$\int\{\hat{\alpha}_{n\ell}(x)-\alpha_{n}(x)\}^{2}F_{0}(dx)=o_{p}(\sigma
_{n}^{2}).$

\bigskip

The next condition imposes rates of convergence for the quadratic remainder%
\[
\hat{\Delta}_{\ell}(w)=\{\hat{\alpha}_{n\ell}(x)-\alpha_{n}(x)\}\{\rho
(w,\hat{\gamma}_{\ell})-\rho(w,\gamma_{0})\}.
\]
Here let $\left\Vert h(W)\right\Vert =\sqrt{E[h(W)^{2}]}.$

\bigskip

\textsc{Assumption A2: }\textit{For each }$\ell=1,...,L$\textit{, either i)}%
\[
\int\hat{\Delta}_{\ell}(w)F_{0}(dw)=o_{p}(\sigma_{n}/\sqrt{n}),\text{ }%
\int\hat{\Delta}_{\ell}(w)^{2}F_{0}(dw)=o_{p}(\sigma_{n}^{2});
\]
\textit{\ or ii) }$\left\Vert \hat{\alpha}_{n\ell}-\alpha_{n0}\right\Vert
\left\Vert \rho(W,\hat{\gamma}_{\ell})-\rho(W,\gamma_{0})\right\Vert
=o_{p}(\sigma_{n}/\sqrt{n}).$

\bigskip

\textsc{Assumption A3: }\textit{i) There is bounded $\alpha_{mn}(X)$ such that}
\[
E[m_{n}(W,\gamma)]=E[v_{mn}(X)\gamma(X)],\text{ \textit{for all} }%
E[\gamma(X)^{2}]<\infty.
\]
ii) \textit{There is }$R(X,\gamma,\gamma_{0})$\textit{ such that}
\[
E[\rho(W,\gamma)|X]=v_{\rho}(X)\{\gamma(X)-\gamma_{0}(X)\}+R(X,\gamma
,\gamma_{0})
\]
\textit{and }$\alpha_{n}(X)=-v_{mn}(X)/v_{\rho}(X)$\textit{;}

\textit{iii) Either }$R(X,\gamma,\gamma_{0})=0$\textit{ or }$E[|\alpha
_{n}(X)R(X,\gamma,\gamma_{0})|]\leq O(\sigma_{n}^{2})\left\Vert \gamma
-\gamma_{0}\right\Vert ^{2},$\textit{ and} $\left\Vert \hat{\gamma}-\gamma
_{0}\right\Vert ^{2}=o_{p}(\sigma_{n}^{-1}/\sqrt{n})$.

\bigskip

\textsc{Lemma A: }\textit{If Assumptions A1 - A3 are satisfied then equation
(\ref{Key}) is satisfied.}

\bigskip

Proof: Let $\phi(w,\gamma,\alpha)=\alpha(x)\rho(w,\gamma),$ $\bar{\phi}%
(\gamma,\alpha)=\int\phi(w,\gamma,\alpha)F_{0}(dw),$ and $\bar{m}_{n}%
(\gamma)=\int m_{n}(w,\gamma)F_{0}(dw).$ Then by adding and subtracting terms
we have
\begin{align*}
&  \frac{1}{n}\sum_{i\in I_{\ell}}\{\psi_{n}(W_{i},\hat{\gamma}_{\ell}%
,\hat{\alpha}_{\ell},\theta_{n})-\psi_{n}(W_{i},\gamma_{0},\alpha_{n}%
,\theta_{n})\}\\
&  =\frac{1}{n}\sum_{i\in I_{\ell}}\{m_{n}(W_{i},\hat{\gamma}_{\ell}%
)+\phi(W_{i},\hat{\gamma}_{\ell},\hat{\alpha}_{\ell})-m_{n}(W_{i},\gamma
_{0})-\phi(W_{i},\gamma_{0},\alpha_{n})\}=\hat{R}_{1}+\hat{R}_{2}+\hat{R}_{3},
\end{align*}
where%
\begin{align}
\hat{R}_{1}  &  =\frac{1}{n}\sum_{i\in I_{\ell}}[m_{n}(W_{i},\hat{\gamma
}_{\ell}-\gamma_{0})-\bar{m}_{n}(\hat{\gamma}_{\ell}-\gamma_{0}%
)]\label{remain}\\
&  +\frac{1}{n}\sum_{i\in I_{\ell}}[\phi(W_{i},\hat{\gamma}_{\ell},\alpha
_{n})-\phi(W_{i},\gamma_{0},\alpha_{n})-\bar{\phi}(\hat{\gamma}_{\ell}%
,\alpha_{n})]\nonumber\\
&  +\frac{1}{n}\sum_{i\in I_{\ell}}[\phi(W_{i},\gamma_{0},\hat{\alpha}_{\ell
})-\phi(W_{i},\gamma_{0},\alpha_{n})],\nonumber\\
\hat{R}_{2}  &  =\frac{1}{n}\sum_{i\in I_{\ell}}[\phi(W_{i},\hat{\gamma}%
_{\ell},\hat{\alpha}_{\ell})-\phi(W_{i},\hat{\gamma}_{\ell},\alpha_{n}%
)-\phi(W_{i},\gamma_{0},\hat{\alpha}_{\ell})+\phi(W_{i},\gamma_{0},\alpha
_{n})]\nonumber\\
&  =\frac{1}{n}\sum_{i\in I_{\ell}}[\hat{\alpha}_{\ell}(X_{i})-\alpha
_{n}(X_{i})][\rho(W_{i},\hat{\gamma}_{\ell})-\rho(W_{i},\gamma_{0}%
)].\nonumber\\
\hat{R}_{3}  &  =\bar{m}_{n}(\hat{\gamma}_{\ell}-\gamma_{0})+\bar{\phi}%
(\hat{\gamma}_{\ell},\alpha_{n}).\nonumber
\end{align}
Define $\hat{\Delta}_{i\ell}=m_{n}(W_{i},\hat{\gamma}_{\ell}-\gamma_{0}%
)-\bar{m}_{n}(\hat{\gamma}_{\ell}-\gamma_{0})$ for $i\in I_{\ell}$ and let
$\mathcal{W}_{\ell}^{c}$ denote the observations $W_{i}$ for $i\notin I_{\ell
}$. Note that $\hat{\gamma}_{\ell}$ depends only on $\mathcal{W}_{\ell}^{c}$
by construction. Then by independence of $\mathcal{W}_{\ell}^{c}$ and
$\{W_{i},i\in I_{\ell}\}$ we have $E[\hat{\Delta}_{i\ell}|\mathcal{W}_{\ell
}^{c}]=0.$ Also by independence of the observations, $E[\hat{\Delta}_{i\ell
}\hat{\Delta}_{j\ell}|\mathcal{W}_{\ell}^{c}]=0$ for $i,j\in I_{\ell}.$
Furthermore, for $i\in I_{\ell}$ $E[\hat{\Delta}_{i\ell}^{2}|\mathcal{W}%
_{\ell}^{c}]\leq\int[m_{n}(w,\hat{\gamma}_{\ell})-m_{n}(w,\gamma_{0}%
)]^{2}F_{0}(dw)$. Then we have
\begin{align*}
E[\left(  \frac{1}{n}\sum_{i\in I_{\ell}}\hat{\Delta}_{i\ell}\right)
^{2}|\mathcal{W}_{\ell}^{c}]  &  =\frac{1}{n^{2}}E[\left(  \sum_{i\in I_{\ell
}}\hat{\Delta}_{i\ell}\right)  ^{2}|\mathcal{W}_{\ell}^{c}]=\frac{1}{n^{2}%
}\sum_{i\in I_{\ell}}E[\hat{\Delta}_{i\ell}^{2}|\mathcal{W}_{\ell}^{c}]\\
&  \leq\frac{1}{n}\int m_{n}(w,\hat{\gamma}_{\ell}-\gamma_{0})^{2}%
F_{0}(dw)=o_{p}(\sigma_{n}^{2}/n).
\end{align*}
The conditional Markov inequality then implies that $\sum_{i\in I_{\ell}}%
\hat{\Delta}_{i\ell}/n=o_{p}(\sigma_{n}/\sqrt{n}).$ The analogous results also
hold for $\hat{\Delta}_{i\ell}=\phi(W,\hat{\gamma}_{\ell},\alpha_{0}%
)-\phi(W,\gamma_{0},\alpha_{0})-\bar{\phi}(\hat{\gamma}_{\ell},\alpha_{0})$
and $\hat{\Delta}_{i\ell}=\phi(W,\gamma_{0},\alpha_{\ell})-\phi(W,\gamma
_{0},\alpha_{0})$ by $\bar{\phi}(\gamma_{0},\hat{\alpha}_{\ell})=0=\bar{\phi
}(\gamma_{0},\alpha_{n}).$ Summing across the three terms in $\hat{R}_{1}$
gives $\hat{R}_{1}=o_{p}(\sigma_{n}/\sqrt{n})$.

Next let $\hat{\Delta}_{\ell}(w)=[\hat{\alpha}_{\ell}(x)-\alpha_{n}%
(x)][\rho(w,\hat{\gamma}_{\ell})-\rho(w,\gamma_{0})].$ Suppose first that
Assumption B (i) is satisfied, so that
\[
\int\hat{\Delta}_{\ell}(w)^{2}F_{0}(dw)=o_{p}(\sigma_{n}^{2}).
\]
Let $\bar{\Delta}_{\ell}=\int\hat{\Delta}_{\ell}(w)F_{0}(dw).$ It then follows
exactly as for $\hat{R}_{1}$ with $\hat{\Delta}_{i\ell}=\hat{\Delta}_{\ell
}(X_{i})-\bar{\Delta}_{\ell}$ that%
\[
\frac{1}{n}\sum_{i\in I_{\ell}}[\hat{\Delta}_{\ell}(X_{i})-\bar{\Delta}_{\ell
}]=o_{p}(\sigma_{n}/\sqrt{n}).
\]
Also by Assumption B i) $\bar{\Delta}_{\ell}=o_{p}(\sigma_{n}/\sqrt{n})$, so
that by the triangle inequality,%
\[
\hat{R}_{2}=\frac{1}{n}\sum_{i\in I_{\ell}}[\hat{\Delta}_{\ell}(X_{i}%
)-\bar{\Delta}_{\ell}]+\frac{n_{\ell}}{n}\bar{\Delta}_{\ell}=o_{p}(\sigma
_{n}/\sqrt{n}),
\]
where $n_{\ell}$ is the number of observations in $I_{\ell}$. Now suppose that
Assumption B ii)\ is satisfied. Then by the triangle and Cauchy-Schwartz
inequalities,%
\begin{align*}
E[\left\vert R_{2}\right\vert |\mathcal{W}_{\ell}^{c}]  &  \leq\int\left\vert
\hat{\Delta}_{\ell}(w)\right\vert F_{0}(dw)\\
&  \leq\sqrt{\int[\hat{\alpha}_{\ell}(x)-\alpha_{n}(x)]^{2}F_{0}(dx)}%
\sqrt{\int[\rho(w,\hat{\gamma}_{\ell})-\rho(w,\gamma_{0})]^{2}F_{0}(dw)}.
\end{align*}
It then follows by the conditional Markov inequality that $\hat{R}_{2}%
=o_{p}(\sigma_{n}/\sqrt{n})$.

Finally, consider $R_{3}$ and let $\hat{\Delta}_{\ell}(x)=\hat{\gamma}_{\ell
}(x)-\gamma_{0}(x)$. By Assumption C,
\begin{align*}
\hat{R}_{3}  &  =\int v_{mn}(x)\hat{\Delta}_{\ell}(x)F_{0}(dx)+\int\alpha
_{n}(x)[\int\rho(w,\hat{\gamma}_{\ell})F(dw|x)]F(dx)\\
&  =\int\{v_{mn}(x)+\alpha_{n}(x)v_{\rho}(x)\}\hat{\Delta}_{\ell}%
(x)F_{0}(dx)\\
&  +\int\left\vert \alpha_{n}(x)R(x,\hat{\gamma}_{\ell},\gamma_{0})\right\vert
F_{0}(dx)\\
&  =0+o_{p}(\sigma_{n}/\sqrt{n})=o_{p}(\sigma_{n}/\sqrt{n}).
\end{align*}
The conclusion now follows by the triangle inequality. $Q.E.D.$

\bigskip

We now use Lemma A to prove Theorems 4 and 6.

\bigskip

Proof Theorem 4: Let $S_{i}=m(W_{i},\gamma_{0})+\alpha_{0}(X_{i}%
)\{Y_{i}-\gamma_{0}(X_{i})\}.$ By the structure of locally linear regression,
as in Fan (1993) equations (2.2)-(2.4), it suffices to prove that%
\begin{align}
\sqrt{nh^{r}}\frac{1}{n}\sum_{i=1}^{n}K_{h}(v-V_{i})[\hat{S}_{i}-S_{i}]  &
=o_{p}(1),\label{rem rate}\\
\sqrt{nh^{r}}\frac{1}{n}\sum_{i=1}^{n}K_{h}(v-V_{i})(v-V_{i})[\hat{S}%
_{i}-S_{i}]  &  =o_{p}(1).\nonumber
\end{align}
For brevity we will provide a proof of only the first expression in
(\ref{rem rate}). The rest of (\ref{rem rate}) can be shown for each of its
elements in an analogous way.

Let
\begin{align*}
m_{n}(W,\gamma)  &  =K_{h}(v-V)m(W,\gamma),\text{ }\alpha_{n}(X)=K_{h}%
(v-V)\alpha_{0}(X),\\
\psi_{n}(W,\gamma,\alpha,\theta)  &  =m_{n}(W,\gamma)+\alpha(X)\{Y-\gamma
(X)\}.
\end{align*}
Let $\sigma_{n}=h^{-r/2}$ and $\hat{\alpha}_{n\ell}=K_{h}(v-V)\hat{\alpha
}(X).$ Note that by $K(u)$ and $\alpha_{0}(X)$ bounded and Assumption 4,
\begin{align*}
\int m_{n}(w,\hat{\gamma}-\gamma_{0})^{2}dF_{0}(w)  &  \leq Ch^{-2r}\left\Vert
\hat{\gamma}-\gamma_{0}\right\Vert ^{2}=o_{p}(\sigma_{n}^{2})\text{,}\\
\int\alpha_{n}(x)^{2}[\hat{\gamma}(x)-\gamma_{0}(x)]^{2}dF(w)  &  \leq
Ch^{-2r}\left\Vert \hat{\gamma}-\gamma_{0}\right\Vert ^{2}=o_{p}(\sigma
_{n}^{2}),\\
\int\{\hat{\alpha}_{n\ell}(x)-\alpha_{n}(x)\}^{2}F_{0}(dx)  &  \leq
Ch^{-2r}\left\Vert \hat{\alpha}-\alpha_{0}\right\Vert ^{2}=o_{p}(\sigma
_{n}^{2}),
\end{align*}
giving Assumption A1.

Next note that $\rho(W,\gamma(X))=Y-\gamma(X)$, so that by Assumption 4,%
\begin{align*}
\left\Vert \hat{\alpha}_{n\ell}-\alpha_{n0}\right\Vert \left\Vert \rho
(W,\hat{\gamma}_{\ell})-\rho(W,\gamma_{0})\right\Vert  &  \leq Ch^{-r}%
\left\Vert \hat{\alpha}_{\ell}-\alpha_{0}\right\Vert \left\Vert \hat{\gamma
}_{\ell}-\gamma_{0}\right\Vert \\
&  =Ch^{-r}o_{p}(h^{r/2}/\sqrt{n})\\
&  =o_{p}(\sigma_{n}/\sqrt{n}).
\end{align*}
so that Assumption A2 is satisfied.

Also Assumption A3 is satisfied with $v_{mn}(X)=K_{h}(v-V)\alpha_{0}(X),$
$v_{\rho}(X)=-1$, and $R(X,\gamma,\gamma_{0})=0$. The conclusion of Theorem 4
then follows from the conclusion of Lemma A with $\sigma_{n}/\sqrt{n}%
=1/\sqrt{nh^{r}}.$ $Q.E.D.$

\bigskip

Proof of Theorem 6: Assumption A1 of Lemma A follows as in the proof of
Theorem 4. Also Assumption A2 is satisfied by Assumption 6 and Assumption A3
by Assumption 7, so the conclusion follows by Lemma A. $Q.E.D.$

\pagebreak

\setlength{\parindent}{-.5cm} \setlength{\parskip}{.1cm}

\begin{center}
\textbf{REFERENCES}
\end{center}

Andrews, I., M. Gentzkow, and J.M. Shapiro (2017): "Measuring the Sensitivity
of Parameter Estimates to Estimation of Moments," \textit{Quarterly Journal of
Economics,} 132, 1553-1592.

Bickel, P.J. (1982): \textquotedblleft On Adaptive
Estimation,\textquotedblright\ \textit{Annals of Statistics} 10, 647--671.

Chaudhuri, P.K. Doksum, and A. Samarov (1997): "On Average Derivative Quantile
Regression," \textit{The Annals of Statistics} 25, 715-744.

Chernozhukov, V., I. Fernandez-Val, and B. Melly (2013): "Inference on
Counterfactual Distributions," \textit{Econometrica }81, 2205-2268.

Chernozhukov, V., D. Chetverikov, M. Demirer, E. Duflo, C. Hansen, W. Newey,
and J. Robins (2018): "Double/Debiased Machine Learning for Treatment and
Structural Parameters," \textit{Econometrics Journal }21, C1--C68.

Chernozhukov, V., W.K. Newey, and R. Singh (2022a): "De-Biased Machine
Learning of Global and Local Parameters Using Regularized Riesz Representers,"
\textit{Econometrics Journal} 25, 576--601.

Chernozhukov, V., W.K. Newey, and R. Singh (2022b): "Automatic Debiased
Machine Learning of Structural and Causal Effects," with V. Chernozhukov and
R. Singh, \textit{Econometrica} 90, 967-1027.

Chernozhukov, V., J.-C. Escanciano, H. Ichimura, W.K. Newey, and J. Robins
(2022c): \textquotedblleft Locally Robust Semiparametric Estimation,"
\textit{Econometrica} 90, 1501-1535.

Chernozhukov, V., W.K. Newey, V. Quintas-Martinas, and V. Syrgkanis (2024):
"Automatic Debiased Machine Learning Via Riesz Regression," https://arxiv.org/pdf/2104.14737.

Conley, T.J., C.B. Hansen, and P.E. Rossi (2012): "Plausibly Exogenous,"
\textit{The Review of Economics and Statistics} 94, 260-272.

Fan, J. (1993): "Local Linear Regression Smoothers and Their Minimax
Efficiencies," The Annals of Statistics 21, 196-216.

Fan, J. and I. Gijbels (1996): \textit{Local Polynomial Modeling and Its
Applications}, Chapman and Hall,\ London.

Firpo, S., N.M. Fortin, and T. Lemieux (2009): "Unconditional Quantile
Regressions,"\textit{\ Econometrica} 77, 953--973.

Foster, D.F. and V. Syrgkanis (2023)\textsc{:} "Orthogonal Statistical
Learning," \textit{Annals of Statistics} 51, 879 - 908.

Gentzkow, M., and J.M. Shapiro (2015): "Measuring the Sensitivity of Parameter
Estimates to Sample Statistics," Working Paper No. 20673, NBER.

Hampel, F. R. (1974): \textquotedblleft The Influence Curve and Its Role In
Robust Estimation,\textquotedblright\ \textit{Journal of the American
Statistical Association} 69, 383--393.

Hausman, J.A. (1978): "Specification Tests in Econometrics,"
\textit{Econometrica} 46, 1251-1271.

Hausman, J. A. and W. K. Newey (2016): \textquotedblleft Individual
Heterogeneity and Average Welfare," \textit{Econometrica} 84, 1225-1248.

Huber, P. (1981): \textit{Robust Statistics,}\ New York: John Wiley and Sons.

Ichimura, H. and W.K. Newey (2022): "The Influence Function of Semiparametric
Estimators," with H. Ichimura, \textit{Quantitative Economics} 13, 29--61.

Imbens, G. and W.K. Newey (2009): "Identification and Estimation of Triangular
Simultaneous Equations Models Without Additivity," \textit{Econometrica} 77, 1481-1512.

Kennedy, E.H. (2023): "Towards Optimal Doubly Robust Estimation of
Heterogeneous Causal Effects," \textit{Electronic Journal of Statistics} 17, 3008-3049.

Kennedy, E.H., S. Balakrishnan, J.M. Robins, L.A. Wasserman (2024): "Minimax
Rates for Heterogenous Causal Effect Estimation," \textit{Annals of Statistics
}52, 793-816.

Khan, S., and D. Nekipelov. "Semiparametric Efficiency in Irregularly
Identified Models." Unpublished working paper (2010).

Koenker, R. and G. Bassett (1978): "Regression Quantiles,"
\textit{Econometrica} 46, 33-50.

Luenberger, D.G. (1969): \textit{Optimization by Vector Space Methods}, New
York: Wiley.

Luedtke, A., and I. Chung (2024): "One-step Estimation of Differentiable
Hilbert-valued Parameters," \textit{Annals of Statistics} 52, 1534-1563.

Machada, J.A.F and J. Mata (2005): "Counterfactual Decomposition of Changes in
Wage Distributions Using Quantile Regression," \textit{Journal of Applied
Econometrics} 20, 445-465.

Melly, B. (2005): "Decomposition of Differences in Distribution Using Quantile
Regression," \textit{Labour Economics} 12, 577-590.

Newey, W.K. (1994): \textquotedblleft The Asymptotic Variance of
Semiparametric Estimators,\textquotedblright\ \textit{Econometrica} 62, 1349--1382.

Newey, W.K. and J.L. Powell (1987): "Asymmetric Least Squares Estimation and
Testing," \textit{Econometrica} 55, 819-847.

Robinson, P.M. (1988): \textquotedblleft Root-N-Consistent Semiparametric
Regression," \textit{Econometrica} 56, 931-954.

Schick, A. (1986): "On Asymptotically Efficient Estimation in Semiparametric
Models," \textit{Annals of Statistics} 14, 1139-1151.

Semenova, V. and V. Chernozhukov (2021):\ "Debiased Machine Learning of
Conditional Average Treatment Effects and Other Causal
Functions\textquotedblright\ \textit{Econometrics Journal} 24, 264-289

Stone, C.J. (1980): "Optimal Rates of Convergence for Nonparametric
Regression," Annals of Statistics, 595-645.

Van der Vaart, A.W. (1991): \textquotedblleft On Differentiable
Functionals,\textquotedblright\ \textit{Annals of Statistics} 19, 178--204.

Van der Vaart, A.W. (1998): \textit{Asymptotic Statistics}, Cambridge:
Cambridge University Press.

Von Mises, R. (1947): \textquotedblleft On the Asymptotic Distribution of
Differentiable Statistical Functions," \textit{Annals of Mathematical
Statistics} 18, 309-348.

\bigskip

\pagebreak

\end{document}